\newcommand{\beq}{\begin{quote}}
\newcommand{\enq}{\end{quote}}
\newcommand{\be}{\begin{equation}}
\newcommand{\en}{\end{equation}}

\newcommand{\De}{\Delta}

\documentstyle[12pt,epsfig]{article}
\begin{document}
\title{ Kepler's Area Law in the
Principia: Filling in some details 
in Newton's proof of Prop. 1.  
} 
\date{}
\author{Michael Nauenberg\\
Department of Physics\\
University of California, Santa Cruz, CA 95064 
}
\maketitle

\begin{abstract}

During the past 25 years there has been a controversy
regarding the adequacy of Newton's proof of
Prop. 1 in Book 1 of the {\it Principia}. This proposition
is of central importance because its proof of
Kepler's area law  allowed Newton to introduce
a geometric measure for time to solve
problems in orbital dynamics in
the {\it Principia}.  It is shown here that the
critics of Prop. 1 have misunderstood Newton's fundamental 
limit argument by neglecting to consider the justification 
for this limit which he gave in Lemma 3.
We clarify the proof of Prop. 1 by filling in some details 
left out by Newton which show that his proof of this proposition 
was adequate and well grounded.

\end{abstract}

{\it Rigor merely sanctions the conquests of sound intuition}, 
Jacques Hadamard


\subsection*{Introduction}

In Prop. 1 of the {\it Principia} Newton gave a proof  that 
Kepler's empirical area law for planetary orbits and the
confinement of these orbits  to a plane are
consequences of his laws of motion for central forces.  
In his words,
\beq
``The areas which bodies made to move in orbits described by
radii drawn to an unmoving center of force lie in unmoving 
planes and are proportional to the times''
\enq
This proposition is justifiably regarded as  
a cornerstone of the Principia, because  
the  proportionality  between  the 
area swept out by the  radius vector of the orbit
and the elapsed  time which this law entails enabled  Newton to solve
dynamical problems by purely  geometrical methods 
supplemented by continuum limit arguments which he had 
developed.  
Although the validity of Newton's proof was
not questioned by  his contemporaries, 
alternative analytic proofs of the area law were given by 
Jacob Hermann and  John Keill  based on the analytic form of the calculus 
which had been developed by  Newton and by  Leibniz (Guicciardini 1999).  During the
past years , however, two influential historians of
science, D.T. Whiteside  and A.J. Aiton
have criticized Newton's proof, claiming  that
it was inadequate and applied only to an
{\it infinitesimal}  arc of the orbit (MP 6:35-36,footnote 19)
(Aiton 1989). Referring to this claim, 
Whiteside remarked  that there were underlying subtleties 
in the proof that Newton  did not fully appreciated, 
and that Newton continued ``to believe in its superficial 
simplicities'' although, Whiteside admitted, not even 
`` Johann Bernoulli, his arch critic ... 
saw fit to impugn the adequacy of Newton's demonstration''.
The basis for the  Aiton-Whiteside criticism  is that
Newton had treated incorrectly the continuum  limit 
of a discrete polygonal orbit due to a sequence 
of central force impulses.  
Subsequently, this criticism  has  been accepted by many
Newtonian scholars  although some  arguments  have  been presented 
that it is  not valid  (Erlichson 1992) (Nauenberg 1998a). 
For example, in his new translation and guide to Newton's Principia, 
I.B. Cohen warmly  endorsed Whiteside's analysis (Cohen 1999)  
while  N. Guicciardini in his new  book {\it Reading the Principia}
questioned whether Newton's limit arguments in Prop.1   
are well grounded, acknowledging, however,  
dissenting views (Guicciardini 1999).
Other recent authors discussing
the {\it Principia}  either neglected to examine 
the validity of Newton's limit arguments
(Brackenridge 1995) (Chandrasekhar 1995), or failed to understand them 
(Densmore 1995).

In this paper we fill in details left out by Newton 
in his presentation of Prop. 1 which we hope will clarify 
the content of this proposition, and help 
resolve the current controversy over Newton's proof.
Specifically, we discuss  his mathematical procedure  
to obtain a continuous orbit as the limit 
of a discrete polygonal orbit, and we show also that in this limit  
a sequence of discrete impulses leads to a continuous 
force with a measure which is proportional to the measure of force  
which Newton gave in  Prop. 6. In this related  proposition, Newton started
directly with the assumption that the orbit is a continuous curve
which satisfies  the area law, and
then obtained a measure for the  central force
by a somewhat different limit argument
from the one which  he presented  in Prop. 1. 
In Prop. 6  this  measure is equal to the acceleration
in units of time proportional to the area swept by the radius
vector.  
In the following discussion the reader should
consult  Props. 1 and 6 which we do not reproduce here  
except for some quotations and the diagrams shown in Figs. 1
and 4. All of these quotations come  from a new English translation 
of the original Latin version of the  {\it Principia} by I.B. Cohen 
and Ann Whitman (Cohen 1999).

The current controversy with Prop. 1 arises partly  because  Newton's 
statement of his limit argument is succinct:

\beq
`` Now let the number of triangles
be increased and their widths decreased indefinitely,
and their ultimate perimeter $ADF$ will (by lemma. 3, corol. 4)
be a curved line ...''
\enq
\begin{figure}
\begin{center}
\epsfxsize=\columnwidth
\epsfbox{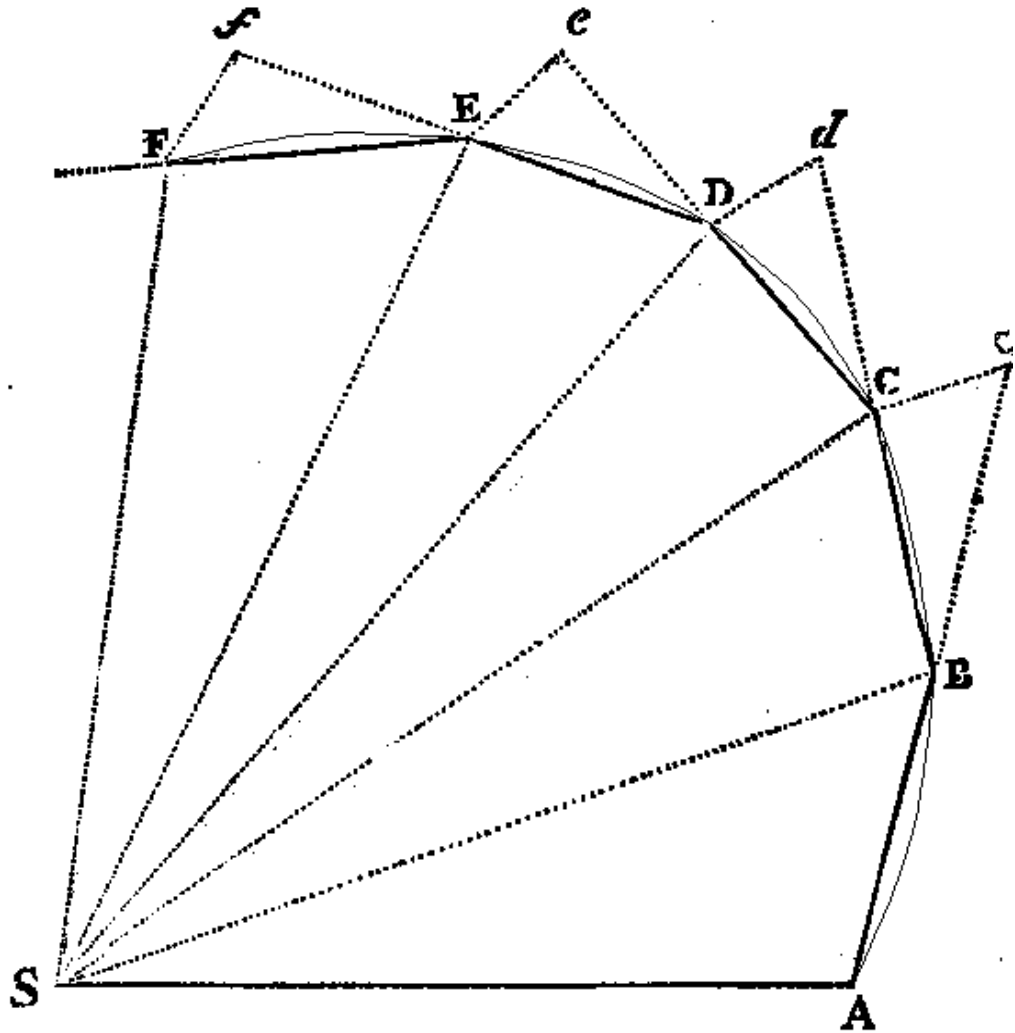}
\end{center}
\caption{
This diagram  corresponds to Newton's diagram  in Prop. 1, 
but with the lines $Sc,Sd,Se$ and $Sf$ deleted, and with 
an additional curve through points $ABCDEF$. 
}
\label{equi3}
\end{figure}

Newton did not give any  details how  this {\it indefinite} increase in
the number of triangles and the corresponding reduction of
their widths would have to be tailored to lead to a well
defined and unique continuum limit.
It is clear  that to understand Newton's  procedure 
one has  to consult Lemma 3 which is given by him  as 
justification for his limit arguments. Remarkably, however, neither  
Whiteside (MP 6:   footnote 19) (Whiteside 1991)   nor
Aiton (Aiton 1989)  commented on this important Lemma,
which was also  neglected  by one of their 
critics (Erlichson 1992), and by other recent  commentators 
of the {\it Principia}  (Brackenridge 1995) (Chandrasekhar 1995) 
(Cohen 1999) (Densmore 1995).
In Lemmas 2 and 3 and its corrolaries, Newton  described 
how the area bounded by a given curve and a line and the
length of the curve   can be approximated 
by a sequence of parallelograms. In Lemma 2, he  
proved rigorously the existence of a limit for this 
area by  obtaining  lower and upper bounds 
given by the area of the inscribed and  circumscribed rectangles
of equal width, showing that for a concave curve
the difference between these two bounds 
is the area of the first rectangle.
Consequently, as the number of these rectangles increases 
indefinitely while their widths approach zero
this difference vanishes, and
the sum of the  area of these rectangles approach the same  limiting
value. By definition, this limit is the area under the 
curve. Indeed, modern calculus books reproduce  Newton's proof for the area under 
a curve, but attribute it to either Cauchy or Riemman. 
In Lemma 3 Newton extended his proof 
in Lemma 2 to the case of rectangles of unequal width,
\beq
``The same ultimate ratios are also ratios of equality when the widths
AB,BC,CD,... of the parallelograms are {\it unequal}[my italics] and are all diminished
indefinitely''.
\enq
It is possible  that this extension 
was included in the {\it Principia} mainly  for its application 
to  Prop. 1, because the  rectangles which
can be associated with the vertices in the corresponding  diagram,
by taking the initial radial position $AS$ as the horizontal axis, see Fig. 1,
would also have unequal widths.
In Cor. 2 of this Lemma Newton asserted that  
\beq
`` ...the  rectilinear figure that is comprehended by the chords
of the vanishing arcs...coincides ultimately with the curvilinear
figure''   
\enq
Finally, in Cor. 4 of Lemma 3 he concluded, 
\beq
``And therefore these ultimate figures (with respect to their
perimeter acE) are not rectilinear, but curvilinear limits of
rectilinear figures.''
\enq 

Hence, Newton's reference to  Lemma 3 is conclusive evidence that 
in Prop. 1 Newton envisioned that 
the vertices of the polygon in his diagram, Fig. 1, 
are located on a given  geometrical curve 
which remains fixed as the number of vertices in this polygon 
increased indefinitely. Newton's construction of this polygon 
requires that these vertices all lie on a plane, and consequently
this curve must  also be  planar.  Unfortunately, this limit curve was not shown 
in the  diagram associated with Prop. 1, and this apparently 
has confused readers who did not consult Lemma 3.  
If any doubt remains for our interpretation, it should be  dispelled 
by Corollaries 2 and 3 of Prop. 1, which were added to
later editions of the {\it Principia}, where Newton called the sides 
$AB$ and $BC$ of the polygon 
\beq

`` chords of two arcs successively described by the same body
in equal times...''. 
\enq
Referring to Fig. 1, it can be seen that with this interpretation  Newton's
entire polygonal construction is  determined by fixing the length
of the first chord $AB$.  
This construction proceeds as follows:   
the extension $Bc$ of this chord
is set equal in length to $AB$, and the first deflection $Cc$ is determined  
by the condition that it is a line  parallel to $BS$ starting at $c$ which intersects
the given curve at the point $C$. This procedure is iterated  
with the next chord $BC$ which is now determined, by
setting  its  extension $Cd$ equal in magnitude to $BC$
and the  deflection $Dd$ parallel to $CS$, starting 
at $d$, and ending  at the intersection $D$ with the given curve. This  
iterative process  continues until the last point $F$ on the curve is reached.
The extension $Bc$ of the chord $AB$ and the deflection $Cc$ must lie
on the  plane of  the initial triangle $SAB$ and therefore the vertex $C$ 
is also on the same plane. Similarly, this property holds also for all 
subsequent vertices of this polygon, or as Newton  stated in his
proposition, 
\beq
``...making the body ...describe  the individual lines $CD,DE,EF,...$,
all these lines will lie in the same plane'' 
\enq
Note, however,  that for these lines to intersect a given curve, 
this  curve  must also be  planar and  lie in this plane. 
The orientation of this plane in space is determined by the initial radial position 
$AS$ and by  the initial chord $AB$ which is in the direction of
the initial velocity for the polygonal orbit. Furthermore, in Newton's 
continuum limit in which  the length of the chord $AB$  vanishes,
the plane's  orientation remains unchanged, 
as the direction of $AB$   
approaches  that of the initial velocity which is directed  
along the tangent of the curve at the initial position at $A$.

In Prop.1 Newton did  not specify directly how the magnitude
of the deflections $Cc, Dd,...$ are obtained, 
nor how the magnitude must vary with the number $n$ of vertices. 
We have seen, however, that
these  deflections can be  determined  by  the assumption that the
polygon vertices are attached to a fixed planar curve.
An alternative possibility, which in fact was considered by Hooke,
(Nauenberg 1994b, 1998b) is that these  deflections 
depend on the radial distance of the vertices of the polygon.
Then to obtain the continuum limit a rule has to be given  for
how these deflections scale with the number of vertices. 
In this case Newton would have had to invoke his curvature
lemma, Lemma 11,  which implies that the deflections
scale as the square of the length of the adjacent chords.
But instead, in Prop. 1  he referred to Lemma 3 which is relevant
to the continuum limit when a fixed curve is given.

Newton's discrete construction in Prop. 1 
refers to a sequence of triangles rather than rectangles  as 
in Lemma 3, but it is easy to see that  this lemma remains also valid in this
case.  There is a detail, however, regarding this application 
which needs some clarification.
As  will be shown below, given a curve of finite length  
Newton's polygonal  construction  does not
in general cover the entire  curve for a finite number of
triangles. This  problem, however, disappears  
in the  limit that the  number of triangles 
increases indefinitely, and the proof is given 
in Appendix A. 

Newton's description of the continuum limit quoted
above continues as follows:

\beq

``...and thus the centripetal force by which the body is continually
drawn back from the tangent of this curve will act uninterruptedly..."
\enq
In  Corols.  3 and 4  to Prop. 1 Newton stated that
the ratio of forces at two distinct point on the curve
where given by the limit of the ratio of the 
displacements caused by central  force impulses at these points, but 
he did not show that in the continuum limit 
the  measure  of this  force is proportional to the  measure for
force which he gave in Prop. 6. This is another detail
that will be discussed after the next section.

\subsection*{The historical context of Prop. 1}

Before we begin our discussion of Prop. 1  it is
necessary to understand what Newton's meant by the term {\it orbit}
when he formulated  Prop. 1,
In the {\it Principia} the term {\it orbit} is not defined explicitly,
but it has been generally understood to mean a geometrical
curve  which describes the position of
a moving body in space. Mathematically  an orbit is a continuous
curve which is  parameterized  by the time variable.
In Prop. 1, however, Newton had to have 
a more restrictive definition, because 
he was dealing there  with the special case 
of the motion of a body under the action of a  central force,
or in his words,

\beq
`` ... bodies made to move in orbits.. by...an unmoving center
of force''.
\enq
To elucidate this point we turn now to the historical circumstances
which led Newton to discover this crucial proposition.
\begin{figure}
\begin{center}
\epsfxsize=\columnwidth
\epsfbox{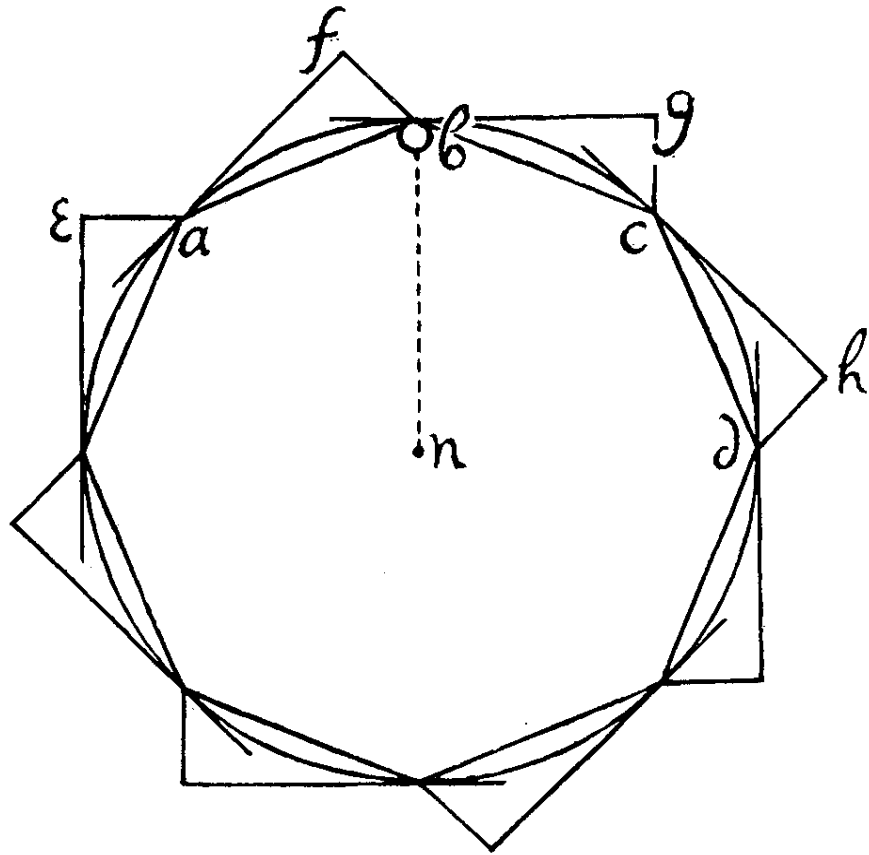}
\end{center}
\caption{
This diagram in Newton's {\it Waste Book} shows
an octagon inscribed in a circle. Here the deflections
$gc,hd$... are shown from the respective tangents
$bg,ch$...rather than from  the extensions of the 
corresponding chords $ab,bc$... as shown  in Fig. 1
}
\label{equi3}
\end{figure}

One of Newton's first ideas about
orbital motion  was to consider  the action of a continuous
force as the limiting case  of a  sequence of force impulses.  
As can be seen from  his earliest surviving  drafts on orbital
motion in the {\it Waste Book} (Herivel 1965) (Whiteside 1991), 
Newton approximated circular motion by a regular polygon 
with its vertices located on a circle, see Fig. 2,.  He also  
obtained an expression for the continuous force 
as the limit of force impulses (Brackenridge 1995).
Apparently, however, he did not generalize this idea to non-circular motion 
until shortly after his correspondence in 1679 with Hooke (Nauenberg 1994b) 
who had suggested in a letter to Newton a somewhat similar 
conceptual scheme to understand
the orbital motion of planets moving around the sun. On Nov. 24, 1679 Hooke  
had written  to Newton 
\beq
`` And particularly if you will let me know your thoughts
of that compounding the celestiall motions of the
planetts of a direct motion by the tangent and 
an attractive motion towards the central body...''
\enq
Indeed, years  later Newton recalled that
\beq
`` In the year 1679 in answer to a letter from Dr. Hook ... I found
{\it now} [my italics] that whatsoever was the law of the forces  which
kept the Planets in their Orbs, the area described by the Radius
drawn from them to the Sun would be proportional to the times
in which they were described...''
\enq
In fact, Prop. 1 appeared for the
first time as Theorem 1  in a short manuscript,  {\it De Motu},
which Newton had sent in 1684 to Halley containing the
beginning draft  of what became later his {\it Principia}.
In this manuscript Newton describe the continuum 
limit in similar words, 
\beq
`` Now let these triangles be infinite in number and infinitely
small, so that each individual triangle corresponds to
the individual moment of time, the centripetal force
acting without diminishing and the proposition will
be established''
\enq
Apart from the mathematical language, which we would regard today 
as far less precise than the language in Prop.1, quoted previously,
it is noteworthy that Newton gave 
no reference here  to any  lemmas 
which would  justify  his limit argument.
Indeed,  he did not give even a hint on how to  
construct  a limiting  procedure in {\it De Motu}. 
Nevertheless, Hooke, who was one of the first members of the Royal
Society to see this manuscript (Nauenberg 1994b,1998b)
recognized  that for a finite number
of impulses Newton's polygonal construction 
gave an approximate solution for orbital motion
along the lines which he had suggested to Newton
in his Nov. 24, 1679 letter quoted above. The
best evidence for this supposition  is that shortly after the
appearance of {\it De Motu},  Hooke implemented 
Newton's construction as an {\it algorithm} to construct
the  orbit when the magnitude of the force impulses is proportional
to the distance from the center.  In  Sept. 1685, 
almost two years before the {\it Principia} was
published, Hooke obtained by this procedure  a remarkably accurate  graphical 
drawing of an  elliptical orbit (Nauenberg 1994b,1998b)
with its center located  at the center of force
by setting the deflections proportional to the distance
from the center. An
enlarged version of the upper part of  his diagram, 
excluding some auxiliary lines, is shown in  Fig. 3. 
\begin{figure}
\begin{center}
\epsfxsize=\columnwidth
\epsfbox{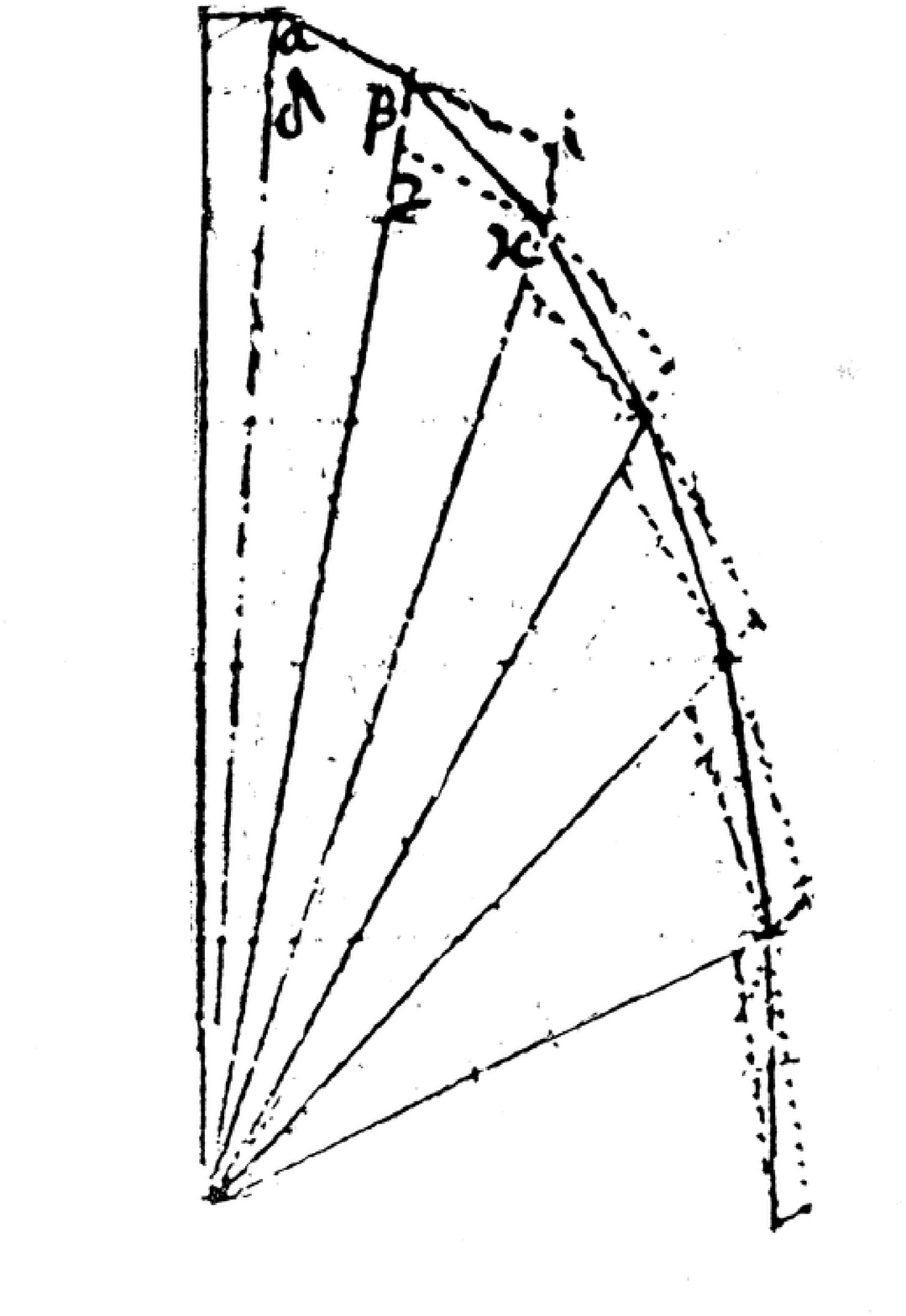}
\end{center}
\caption{
This diagram is a blowup of the upper
part of Hooke's Sept. 1685  diagram (see reference \cite{hooke1}
with some auxiliary lines deleted to show 
more clearly its correspondence with Newton's diagram
in Prop. 1, see Fig. 1. 
}
\label{equi3}
\end{figure}

This reveals quite clearly the relation of Hooke's diagram  
to Newton's  diagram in Prop. 1, which is shown in Fig. 1
The main  difference is that in Hooke's figure the initial 
position is located above the center of force and the initial  velocity 
points to the rigth,  leading to clockwise motion, while
in Newton's figure the initial position is to the rigth
of the center of force and the initial  velocity is mainly 
upwards leading to counter-clockwise motion.
Indeed,  Hooke had  also conjectured that the force
of gravity consisted of discrete pulses. 
In one of his Cutlerian  lectures entitled a 
{\it Discourse on the Nature of Comets}, read at a meeting of
the Royal Society soon after Michaelmas 1682, but
published only after his death,
Hooke speculated that bodies emitted periodic gravitational
pulse in analogy with his theory of sound and light,
and deduced that the intensity decreased with the inverse
square distance from the source:
\beq
``This propagated Pulse I take to be the Cause of the 
Descent of Bodies towards the Earth... Suppose  for Instance
there should be  a 1000 of these pulses in a Second of
Time, then must the Grave body receive all those thousand 
impressions within the space of that Second, and a thousand
more the next...(Hooke 1705)''
\enq 
But the important  question of how the magnitude of the
deflection caused by the force  impulses 
scales with the size of the triangles, which is 
essential to establish the existence of  a continuum limit
in this application of Newton's polygonal construction,  
was not - and undoubtedly could not be -  raised by Hooke. 
Only later did Newton show, in his Lemma 11 on curvature which appeared
in  section 1 of the {\it Principia}, 
how these deflections scale with the size of the adjacent 
chord or arc length, namely as the square of these quantities.

From the foregoing it is therefore  reasonable to conclude  
that in Prop. 1 Newton had in mind that any orbit
under the action of central forces is the continuum  limit of
a polygonal orbit, which is {\it caused} by  a sequence of force impulses. 
In this case the  body moves along straight lines  between impulses, 
or by ``direct motion'' as envisioned by  Hooke, where this
motion is along  the sides of a polygon, 
while the force impulse give rise  to a linear deflection 
or an ``attractive motion'' towards the center. 
Since all the impulses  are directed to this  common center the  resulting  
polygonal orbit must be in a plane as Newton demonstrated in 
Prop. 1. The orientation of this
plane is determined by the direction of the  velocity and 
the position of the body relative to the center of force at some initial time. 
Newton was well aware of the important role of  initial conditions
to fix the orbit, and in Prop.17 he discussed these conditions
for the case of elliptical motion under the action of
inverse square forces,

\beq
`` Suppose that the centripetal force is inversely proportional
to the square of the distance of places from the center ...
it is required to find the  line which a body describes when
going forth from a given place with a given velocity along
a given straight line.''
\enq

Setting  the  sequence of  central force impulses 
at {\it equal} time intervals,  Newton gave a proof in Prop. 1  
that the  areas of the triangles associated with the  
resulting polygonal orbit are equal.
Since  planarity as well as the 
this  area law  are properties of any polygonal orbit due
to central force impulses, it is reasonable to expect
that these properties remain also  valid 
in a {\it properly} defined  continuum limit, but Newton 
did not give any details about how this limit is obtained apart from
the brief sentence quoted in our Introduction, and his reference to Cor. 4 Lemma 3. 
In the next section  we will attempt to fill in some of the details left out 
in Newton's discussion. We  should keep in mind 
that for any finite polygon  the time associated with the vertices
of the polygon is different from the corresponding time 
associated with  the continuous orbit. Newton was not always careful
about this distinction. For example, 
in Cor. 2 of Prop.1 quoted above, Newton considered   
``two arcs successively  described in equal times'', when 
he was  referring  to  the time intervals  
defined by his polygonal construction,
but for a finite polygon these time intervals differ.   
The reason is that  the area of a finite triangle  differs by a small amount
from the area of the corresponding ``pie'' which is  bounded by an arc  
instead of a chord of the orbital curve, and
these  areas are  proportional to time intervals.

\subsection*{ Filling in some details of Newton's proof of Prop. 1}

Referring to the diagram in Prop. 1, see Fig.1, we assume that the
vertices  $A$, $B$, $C$, $D$, $E$ and $F$ of Newton's  polygon 
are  located on a given  curve. Since this polygon is planar
we expect  that this curve should also lie on the same plane
We  will show
that in the continuum limit, Newton's polygonal construction  
determines a parameterization of this curve as a function of time, 
describing  orbital motion under the action of a central
force centered at the point $S$, and that in this limit  
the magnitude  of the central force obtained
from Prop. 1 is equivalent  to that defined in Prop. 6. 
There are some restrictions on the possible planar curves which
can support Newton's polygonal construction. For example, the radius vector 
$\vec r$ with origin at $S$ cannot become tangential 
to the curve, because in the neighborhood of any
such point  Newton's  polygonal approximation cannot be 
constructed. This construction
also fails  when the curve
crosses this  origin, which  corresponds to
orbital motion  when the central force diverges  
as $1/r^3$ or faster, and when the curvature approaches
infinity. Therefore, our discussion  
will be confined to regions of space where the  
central force and the
curvature of the orbit remain finite.

As shown in the Introduction, given the length of the initial 
chord $AB$ it is evident that 
Newton's polygonal construction is uniquely  specified
by the condition that the  vertices of the polygon
lie on a {\it fixed} planar curve.  That Newton had such a  given curve 
in mind, although it did not appear
in the diagram in Prop. 1, is clear from his reference 
to Corol. 4, Lemma 3 for the continuum limit, as we
argued in detail previously. 
While in this lemma the approximation to  a continuous curve  
is discussed for a subdivision in rectangles, the extension 
to triangles is quite  straightforward, but
there is a detail which needs to be worked out:
if $A$ is the initial point of the curve then  
unless the length of the initial chord $AB$ is suitably chosen
the last point $F$ of the polygon  will in general
not lie at the  endpoint of the given  curve. 
But in the continuum limit  this is not a problem. 
Suppose that the  last vertex $F$
occurs before the endpoint of the curve.  Apply  Newton's construction
by extending  chord $EF$ to $g$,  and draw a line from $g$ parallel
to $FS$. Then either a) this line intersects the
curve at a new vertex $G$ or b) it does not intersect the
curve at all. In case a) repeat Newton's  construction until
case b) is reached. When  $F$ is the last vertex,
and deflections due to the central force impulse are small,
it is expected that the  distance of $F$ from the endpoint of the curve
decreases as the length of the initial chord $AB$ is decreased. 
Then in the continuum limit all  the cord lengths becomes vanishingly small,
and the last point $F$ reaches the end point of the curve.
A rigorous proof for this assertion is given in Appendix A  
which is based on a suggestion by Bruce Pourciau (private communication). 
These considerations are valid provided that
the curvature of the orbit is finite, in which case the difference
between the chord and the arc length is second order in the chord length. 

In Prop. 1, the deflection  at each vertex 
due to the central force impulse    
is similar to the deflection that Newton described in Prop. 6. 
The main difference is that in Prop. 1
the extension at the end of the adjacent chord 
replaces the tangent to the curve in Prop. 6. For example, 
referring to Figs. 1 and 4,  at vertex  $B$ the extension $Bc$ of 
the chord $AB$ in Prop. 1
corresponds to the tangent line $PR$ in Prop. 6, with $P$ equivalent
to $B$, and the deflection
$Cc$ which is parallel to $BS$ in Prop. 1 corresponding  
to  the deflection $RQ$ 
which is parallel to $PS$ in Prop. 6.
In the limit that the chord length approaches zero, the
difference between the tangent and the chord becomes vanishingly small,
and consequently these two  constructions become similar,
except that the magnitude of the deflection at a vertex in Prop. 1 is 
twice as large as that in Prop. 6.
Hence, in the continuum limit  Newton would have
obtained a  measure for central force in Prop. 1 
equivalent to that in Prop. 6, 
by dividing the deflection at a vertex, which depends  quadratically
on the adjacent chord length, by the square of twice the area 
of the triangles.  
For example,  the measure  of the force at $B$ which is 
the continuum limit of  $Cc/(\De SBC)^2$,
where $\De SBC$ is twice the area of the triangle $SBC$, 
is twice the measure  of the force at P  
given in Prop. 6, which is the limit of $QR/(\De SPR)^2$, 
where $\De SPR =SP \times QT$
is twice the area of the triangle $SPR$. 
\begin{figure}
\begin{center}
\epsfxsize=\columnwidth
\epsfbox{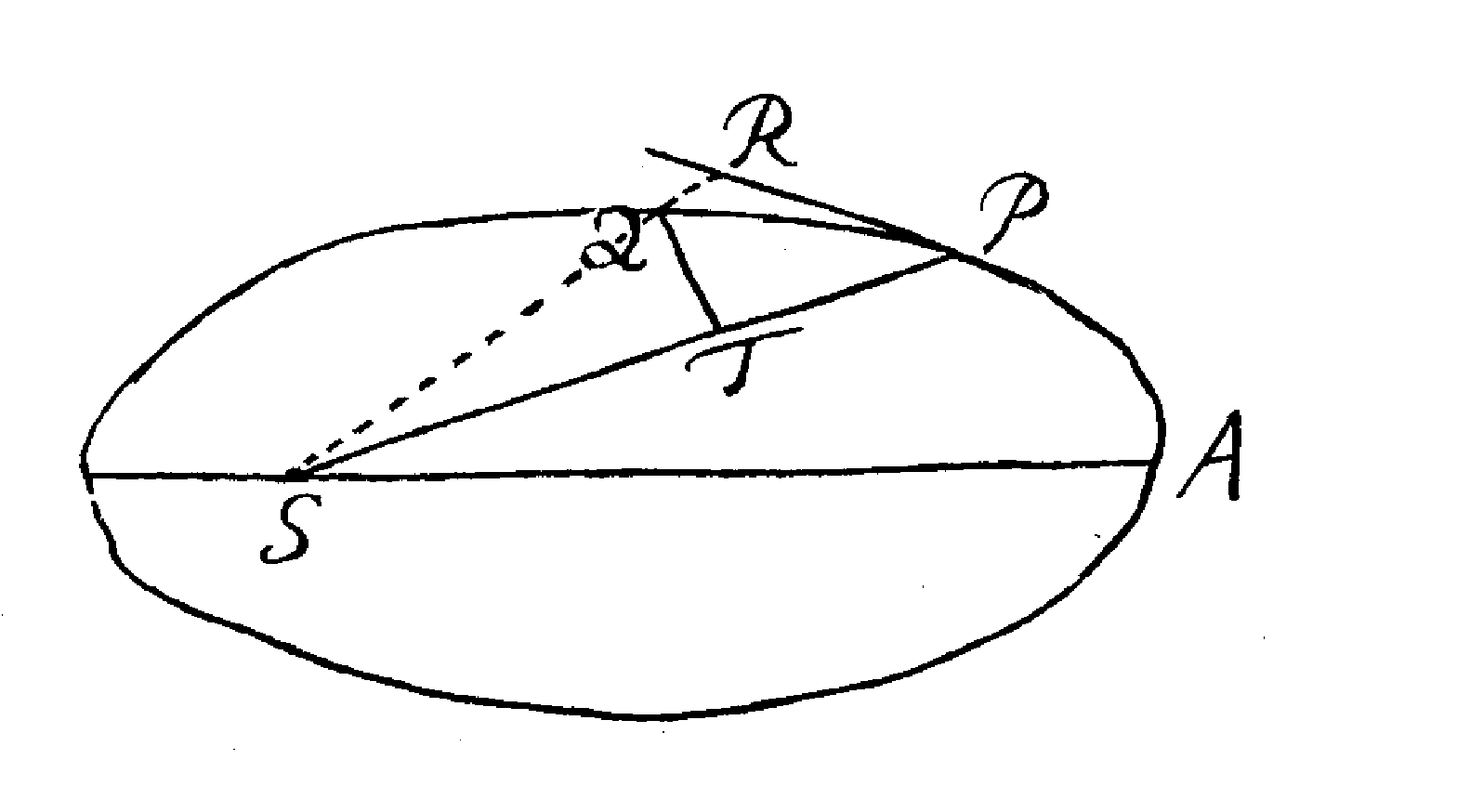}
\end{center}
\caption{
Newton's diagram for Prop.6 
}
\label{equi3}
\end{figure}

\subsection* {Conclusion}
We have shown that apart from some details which we
have discussed here,  Newton's polygonal
construction in Prop. 1, Fig. 1, has a well defined
continuum  limit which is  justified
by Lemmas 2 and 3 in section 1 of the {\it Principia}. This limit
parameterizes a continuous planar curve as a function of 
time which describes   orbital motion under the
action of central forces with origin at $S$.
This orbit  satisfies a generalization of  Kepler's  area law 
which states that the time interval between any two points on the orbit 
is proportional to the area swept out by the corresponding radius vector.
This definition of an orbit for central forces is the one  
which Newton introduced in Prop. 6 to evaluate the 
magnitude of the central force. 
The planarity property  of the orbit is a straightforward  consequence
of the requirement that an  orbital curve is the limit
of a Newtonian polygonal construction,
because the vertices of this polygon  all
lie in the same plane  when the force impulses are directed to a common center.
The orientation of this plane is determined entirely by some initial 
conditions, i.e., 
the position and velocity of the moving body at some given time.
We also have  shown that in the continuum limit
these  impulses lead to a continuum central force 
which is proportional  to the force measure  described in Prop. 6.

Historically, Hooke played an important role  in prompting Newton in 1679 
to take  a new  approach to orbital dynamics which led  him to prove 
the area law for central forces (Nauenberg 1994b) which Kepler had
found  empirically (by fitting the orbit of the
planet Mars to the observations of Tycho Brahe).
There is evidence that by 1679  Newton had been pursuing a different
approach to orbital dynamics based on his development
of curvature (Nauenberg 1994a, Brackenridge and Nauenberg 2001).
This led him to a local description of
central forces in which the area law was not apparent.
We conclude that apart from some details left out by Newton, 
which have been discussed here, Prop. 1  
is well grounded and provides a valid proof
for central force of the generalization of  Kepler's area law 
and the  property  that the orbits lie in a plane. 

\subsection*{Appendix A}

Following a suggestion by Bruce Pourciau (private communication),  
we give here a rigorous proof that for an orbit with 
finite curvature all the chord lengths in Newton's polygonal
construction in Prop. 1 vanish in a mathematically
well defined  continuum limit, and that in this limit 
this construction covers a given length of this orbit.
This condition is neccesary for the application of Lemma 3
to Prop. 1.

Let $s(j)$ be the cord length associated with the
jth and (j+1)th vertices,  and $e(j)=s(j+1)-s(j)$
the difference in length between adjacent chords.
Then $e(j)\approx d(j)cos(\theta(j)$ where $d(j)$ is the
magnitude of the deflection at the jth vertex
and $\theta (j)$ is the angle between this deflection
and the adjacent chord. For an orbit with
finite curvature $d(j)\approx  c(j) s(j)^2/sin(\theta (j))$, 
where the constants $c(j)$ are
proportional to the curvature at the j-th vertex as $s(j)$ approaches zero.
Given  an orbital  curve,  the length of the first chord $s(1$ 
then determines the length of  all the other cords in  
Newton's polygonal construction, where 

\be
\label{chord1}
s(j)=s(1)+e(1)+e(2)+...+e(j-1)
\en
for $j=2,3,...n$.
Let  $s1=L/n$ where $L$ is some fixed length. Then
\be
e(1)=c'(1)*(L/n)^2,
\en
and for $j=2,3,...n$
\be
e(j)=c'(j)*(L/n)^2 + o(l/n)^3
\en
where $c'(j)=c(j)*cot(\theta (j))$ and
$o(L/n)^3$ refer to all terms proportional to $(L/n)^3$ and higher powers
of $(L/n)$.
Hence, the sum 
\be
e(1)+e(2)+...+e(j)= (L/n)^2*(c(1)+c(2)+...+c(j)) + o(L/n)^3
\en

Now suppose that c is the  maximum curvature of the orbit.
Then 

\be
e(1)+e(2)+...+e(j)<(L/n)^2*(j-1)*c + o(L/n)^3, 
\en
and
\be
e(1)+e(2)+...+e(n)<(L/n)^2*(n-1)*c+ o(L/n)^2  
\en

As  $n$ goes to infinity both of these two sums go to zero,
and therefore, according to  Eq. \ref{chord1}  
all the chords $s(j)$ vanish in this limit. Therefore in this limit
Newton's polygonal construction in Prop. 1covers any chosen  segment of the orbit  
with a length given by the limit of the sum $s(1)+s(2)+s(3)+...+ s(n)$.
This length  depends on the parameter $L$, 
but it is equal to $L$ only in the special case of a circular
orbit. In this case the angle $\theta(j)$ is equal $90^o$, and the $e(j)'s$
vanishes as  $(L/n)^3$

\subsection* {Appendix B}
We expresses the content of Prop.1 in modern vector notation,
which may be  helpful in clarifying some of the main points which we emphasized 
previously.
Assuming that there are $n$ vertices in the polygon in Fig. 1 of Prop. 1,
let $\vec r_j$ be the position vector and $\vec v_j$ the velocity  
vector at the j-th vertex
where $j=1,2,...n$. Then Newton's construction take the  the form
\be
\label{pos1}
\vec r_{j+1}=\vec r_j+\vec v_j \De t,
\en
and 
\be
\label{vel1}
\vec v_{j+1}=\vec v_j+\vec \De v_{j+1},
\en
where $\De t$ is the {\it equal} time interval between impulses  introduced by
Newton, and $\De \vec v_j$ is instantaneous velocity change 
$\De \vec  v_j= \vec d_j/\De t$ where $\vec d_j$ 
is the deflection at the
j-th vertex. Then the crossed product of Eqs. \ref{pos1} and \ref{vel1}
\be
\label{area1}
\vec r_{j+1} \times \vec v_{j+1}=
\vec r_j \times \vec v_j +\vec r_{j+1} \times \De \vec v_{j+1}
\en
where $(1/2)|\vec r_j \times \vec v_j| $ is the area of the  triangle associated
with the j-the vertex. Since  the deflection $\vec d_j $ and corresponding 
velocity change $\De \vec v_j$ are parallel to $\vec r_j$,
which is Newton's {\it definition} of central force impulses in Prop. 1, 
then the last term in Eq. \ref {area1} vanishes, and consequently
a) the area of the triangles are equal, and b) the vertices of the polygon lie 
on a plane. This constitutes  Newton's proof of Prop. 1 in the language
of  vector calculus.

The problem of the  continuum limit is  to describe  
how the time interval  $\De t$ and the deflections $\vec d_j$  
should  vary as $n$ approaches  infinity.
Newton's   statement 
\beq
``Let the time be divided in equal times...'',
\enq
corresponds  to setting  $\De t=T/n$, where $T$ is some finite time interval,
and his reference to  Cor. 4 Lemma 3 implies that the deflections 
$\vec d_j$ are determined by the 
condition that the vertices of the polygon are located on a {\it given}
orbital curve. This  curve can be described  by a vector $\vec R(u)$ where $u$ is 
a scalar parameter which  can be chosen arbitrarily. For example, $u$ can be the area swept by
this radius vector with origin at the center of force in which case it 
will turn out to be  proportional
to time for the continuum orbit. But it can also be the arc length, 
or the angular variable in polar coordinates which are  non-linear functions of time.
Then the condition  that the j-th vertex is located on this curve is given by 
\be 
\vec r_j=\vec R(u_j),
\en
where $u_j$ is the value of $u$ at the position of this vertex.
It can be readily seen that in this case
\be
\vec d_j =\vec R(u_{j+1})+\vec R(u_{j-1})-2 \vec R_j,
\en
and as $\De t$ vanishes the ratio 
\be
\frac {d_j}{(v_j \De t)^2}
\en
where $d_j=|\vec d_j|$ and $v_j=|\vec v_j|$,
is  proportional to the curvature of the orbit at time t. 
Hence this ratio as a well defined limit
in the case that this curvature is finite. 
Likewise,
\be
\frac {\vec \De v_j}{\De t} =\frac{\vec d_j}{\De t ^2}.
\en
has a continuum limit corresponding to the
standard calculus definition of the acceleration or force/unit mass  $\vec a$
\be
\vec a=\lim_{n\to\infty} \frac{\vec \De v_j}{\De t}=\frac {d^2 \vec R}{dt^2} 
\en

\subsection*{Acknowledgements}

I would like to thank Bruce Pourciau for many stimulating
exchanges, suggestions and critical comments on the subject of this paper.

\end{document}